\documentclass[11pt,reqno,a4paper]{amsart}
\usepackage{amsmath}
\usepackage{amssymb}
\usepackage{enumerate}
\newcommand{\eps}{\varepsilon}
\newcommand{\N}{\mathbb{N}}
\newcommand{\Z}{\mathbb{Z}}
\newcommand{\Q}{\mathbb{Q}}
\newcommand{\R}{\mathbb{R}}
\newcommand{\T}{\mathbb{T}}

\newcommand{\pf}{f_*\mu}

\newcommand{\pep}{\psi_{x,\varepsilon_n}(y)}
\newcommand{\pdep}{\psi_{x,2\varepsilon_n}(T^ny)}
\newcommand{\pdepy}{\psi_{x,2\varepsilon_n}(y)}
\DeclareMathOperator{\diam}{diam}
\DeclareMathOperator{\rank}{rank}
\renewcommand{\P}{\mathbb{P}}

\title{Poincar\'e Recurrence for observations}
\author{J\'er\^ome Rousseau and Beno\^\i t Saussol}
\address{Universit\'e Europ\'eenne de Bretagne, Universit\'e de Brest, Laboratoire de Math\'ematiques UMR CNRS 6205, 6 avenue Victor le Gorgeu, CS93837, F-29238 Brest Cedex 3}
\email{jerome.rousseau@univ-brest.fr, benoit.saussol@univ-brest.fr}
\urladdr{http://www.math.univ-brest.fr/perso/benoit.saussol}
\date{July 4, 2008}

\begin{document}
\newtheorem{lem}{Lemma}
\newtheorem{theorem}[lem]{Theorem}
\newtheorem{defi}[lem]{Definition}
\newtheorem{prop}[lem]{Proposition}
\newtheorem{corollary}[lem]{Corollary}
\newtheorem{example}[lem]{Example}
\newtheorem*{remark}{Remark}

\maketitle

\begin{abstract}
A high dimensional dynamical system is often studied by experimentalists through the measurement of a relatively low number of different quantities, called an observation. Following this idea and in the
continuity of Boshernitzan's work, for a measure preserving
system, we study Poincar\'e recurrence for the observation.
The link between the return time for the observation and the Hausdorff dimension of the
image of the invariant measure is considered. We prove that
when the decay of correlations is super polynomial, the recurrence rates for the observations and the pointwise dimensions relatively to the push-forward are equal.
\end{abstract}

\section{Introduction}
The famous Zermelo paradox reveals that the classical Poincar\'e recurrence theorem has some implications out of physical sense. Indeed, if we start with all the particles in one side of a box, nobody will ever see all the particles coming back in one side of our box at the same time. Nevertheless, if we focus on a few number of these particles, this event will appear after a reasonable time. In the same way, when we study a high dimensional dynamical system we might not know all the aspects of the evolution but only a part or certain quantities of the system. This might be due to the difficulty to study a high dimensional system, but also to the lack of interest of an over-detailed description.

Recently, Ott and York tried to elaborate some Platonic formalism of dynamical systems \cite{MR2031277}. The reality, the dynamical system $(X,T,\mu)$, is only known through a measurement or observation, that is a function defined on $X$ taking values in (typically) a lower dimensional space. The following result by Boshernitzan \cite{MR1231839} about Poincar\'e recurrence falls in this frame. If we have a measure preserving dynamical system $(X,T,\mu)$ and an observable $f$ from $X$ to a metric space $(Y,d)$ then whenever the $\alpha$-dimensional Hausdorff measure is $\sigma$-finite on $Y$ we have
\begin{equation}\label{bosh}\liminf_{n\rightarrow\infty}n^{1/\alpha}d\left(f(x),f(T^nx)\right)<\infty\qquad\textrm{for $\mu$-almost every $x$}.\end{equation}

The main aim of this paper is to prove a refinement of \eqref{bosh} and a generalization of \cite{MR1833809,MR2191396} for recurrence rates for observations.

In Section \ref{statement}, we give the precise definition of the recurrence rates for the observations and state an upper bound in term of dimension (Theorem~\ref{thsup} which is proved in Section~\ref{majoration}), then under an additional assumption we state our main result (Theorem~\ref{thprinc2} which is proved in Section~\ref{egalite}), and finally, we analyze in the case of the Lebesgue measure the existence of the pointwise dimension for its smooth image (Theorem~\ref{thfcinf} which is proved in Section~\ref{dimensions}).
\section{Statement of the results}\label{statement}
\subsection{Definitions and general inequality}
Let $(X,\mathcal{A},\mu,T)$ be a measure preserving system (m.p.s.) i.e. $\mathcal{A}$ is a $\sigma$-algebra, $\mu$ is a measure on $(X,\mathcal{A})$ with $\mu(X)=1$ and $\mu$ is invariant by $T$ (i.e $\mu(T^{-1}A)=\mu(A)$ for all $A\in\mathcal{A}$) where $T:X\rightarrow X$.\\
Let $f:X\rightarrow Y$ be a function, called observable (we will specify the space $X$ and $Y$ later). We introduce the return time for the observation and its associated recurrence rates.
\begin{defi}Let $f:X\rightarrow Y$ be a measurable function, we define for $x\in X$ the return time for the observation:
\[\tau_r^f(x):=\inf\left\{k\in\N^*,\,f(T^kx)\in B\left(f(x),r\right)\right\}\]
where $B(x,r)$ is the ball centered in $x$ with radius $r$. We then define the lower and upper recurrence rate for the observation:
\[\underline{R}^f_i(x):=\liminf_{r\rightarrow0}\frac{\log\tau_r^f(x)}{-\log r}\qquad\overline{R}^f_i(x):=\limsup_{r\rightarrow0}\frac{\log\tau_r^f(x)}{-\log r}.\]
We also define for $p\in\N$ the $p$-non-instantaneous return time for the observation:
\[\tau_{r,p}^f(x):=\inf\left\{k>p,\,f(T^kx)\in B\left(f(x),r\right)\right\}.\]
Then we define the non-instantaneous lower and upper recurrence rates for the observation:
\[\underline{R}^f(x):=\lim_{p\rightarrow\infty}\liminf_{r\rightarrow0}\frac{\log\tau_{r,p}^f(x)}{-\log r}\qquad\overline{R}^f(x):=\lim_{p\rightarrow\infty}\limsup_{r\rightarrow0}\frac{\log\tau_{r,p}^f(x)}{-\log r}.\]
Whenever $\underline{R}^f(x)=\overline{R}^f(x)$ we denote by $R^f(x)$ the value of the limit.
\end{defi}
\emph{The lower and upper pointwise or local dimension} of a Borel probability measure $\nu$ on Y at a point $y\in Y$ are defined by
\[\underline{d}_\nu(y)=\underset{r\rightarrow0}{\underline\lim}\frac{\log\nu\left(B\left(y,r\right)\right)}{\log r}\qquad\textrm{and}\qquad\overline{d}_\nu(y)=\underset{r\rightarrow0}{\overline\lim}\frac{\log\nu\left(B\left(y,r\right)\right)}{\log r}.\]
The pushforward measure $f_*\mu(.):=\mu(f^{-1}(.))$ is a probability measure on $Y$ and we define \emph{the lower and upper pointwise dimension for the observations} with respect to $\mu$ at a point $x\in X$ by
\[\underline{d}^f_\mu(x)=\underline{d}_{\pf}(f(x))\qquad\textrm{and}\qquad\overline{d}^f_\mu(x)=\overline{d}_{\pf}(f(x)).\]
If they are equal, we denote by $d_\mu^f(x)$ the common value.
\begin{theorem}\label{thsup}
Let $(X,\mathcal{A},\mu,T)$ be a m.p.s, let $f:X\rightarrow \R^N$ be a measurable function. Then
\[\underline{R}^f(x)\leq\underline{d}^f_\mu(x)\qquad\textrm{and}\qquad\overline{R}^f(x)\leq\overline{d}^f_\mu(x)\]
for $\mu$-almost every $x\in X$.
\end{theorem}
This result is satisfactory in the sense that it holds for any dynamical system and observation. Moreover, under natural assumptions we will show that the equality is true. Still, these inequalities may be strict, the caricatural example is when $T$ is the identity map.
\begin{example}Let $(\Omega,\mathcal{F}, \P)$ be a probability space together with a $\P$-preserving map $\theta$ and $Y\subset\R^N$ a Borel set. The family $(F_\omega)_{\omega\in\Omega}$ is called a random transformation, where for each $\omega$, $F_\omega$ is a map from $Y$ to $Y$ such that the map $(w,y)\rightarrow F_\omega(y)$ is $\mathcal{F\times B}(\R^N)$-measurable. The map $T:X=\Omega\times Y\rightarrow X$ defined by $T(\omega,y)=(\theta\omega,F_\omega(y))$ is called a skew product transformation. Let $\mathcal{M}_\P(X,T)$ be the set of $T$-invariant probability measure having the marginal $\P$ on $\Omega$. For any $\mu\in\mathcal{M}_\P(X,T)$, Theorem~\ref{thsup} applies with $f$ the projection on $Y$, and gives an upper bound for the time needed by a typical random orbit $F_{\theta^k\omega}\circ...\circ F_{\theta\omega}\circ F_\omega(y)$ to come back  close to its starting point $y$.
\end{example}
\subsection{Poincar\'e recurrence for observations}
{}From now on, let assume that $X$ is a metric space and $\mathcal{A}$ is its Borel $\sigma$-algebra. We can then introduce the decay of correlations:
\begin{defi}
$(X,T,\mu)$ has a super-polynomial decay of correlations if, for all $\phi$, $\psi$ Lipschitz functions from $X$ to $\R$ and for all $n\in\N^*$,  we have:
\[\left|\int_X\phi\circ T^n\,\psi d\mu-\int_X \phi d\mu\int_X\psi d \mu\right|\leq\|\phi\|\|\psi\|\theta_n\]
with $\lim_{n\rightarrow\infty}\theta_n n^p=0$ for all $p>0$ and where $\|.\|$ is the Lipschitz norm.
\end{defi}
The main result of our paper is:
\begin{theorem}\label{thprinc2}
Let $(X,\mathcal{A},\mu,T)$ be a m.p.s with a super-polynomial decay of correlations. Let $f:X\rightarrow \R^N$ be a Lipschitz observable. Then, we have
\[\underline{R}^f(x)=\underline{d}^f_\mu(x)\qquad\textrm{and}\qquad\overline{R}^f(x)=\overline{d}^f_\mu(x)\]
for $\mu$-almost every $x$ such that $\underline{d}^f_\mu(x)>0$.
\end{theorem}
Taking the identity function for $f$, we recover the result of \cite{MR1833809} and \cite{MR2191396} under weaker assumptions. The main assumption of the theorem about decay of correlations is satisfied in a variety of systems with some hyperbolic behavior and studied in an abundant literature (e.g. \cite{MR1637655,MR1665752,MR1793194}).
\begin{defi}We say that a probability measure $\nu$ is exact dimensional if there exists a constant $d_\nu\in\R$ such that
\[\underline{d}_\nu(\cdot)=\overline{d}_\nu(\cdot)=d_\nu\qquad\textrm{almost everywhere.}\]
\end{defi}
It is well known that in this case many notion of dimension coincide (see Section~\ref{sechaus} for details). In particular the Hausdorff dimension $\dim_H\nu$ satisfies
\begin{prop}\label{exdim} If $\nu$ is exact dimensional, then
\[d_\nu(\cdot)=\dim_H\nu\,\,\textrm{ almost everywhere}.\]
\end{prop}
\begin{corollary}\label{thprinc}
Let $(X,\mathcal{A},\mu,T)$ be a m.p.s with a super-polynomial decay of correlations. Let $f:X\rightarrow \R^N$ be a Lipschitz observable. Then, if $f_*\mu$ is exact dimensional, we have
\[\underline{R}^f(x)=\overline{R}^f(x)=\dim_H\pf\qquad \textrm{for $\mu$-almost every $x\in X$.}\]
\end{corollary}
\begin{remark}We have the equivalence
\[f_*\mu\textrm{ is exact dimensional }\Longleftrightarrow\exists d,\,\underline{d}^f_\mu(x)=\overline {d}^f_\mu(x)=d\qquad\textrm{for $\mu$-almost every $x\in X$}.\]
\end{remark}
\begin{proof}[Proof of Corollary~\ref{thprinc}]If $\dim_H\pf=0$, then the conclusion follows from Theorem~\ref{thsup} and Proposition~\ref{exdim}. In the general case, it is just a combination of Theorem~\ref{thprinc2} and Proposition~\ref{exdim}. \end{proof}
Theorem~\ref{thprinc2} does not apply to those points where $\underline{d}^f_\mu(x)=0$. When $\overline{d}^f_\mu(x)=0$ also, this is not a restriction because Theorem~\ref{thsup} applies and gives $\overline{R}^f(x)=\underline{R}^f(x)=0$. However, the question remains when $\overline{d}^f_\mu(x)\neq\underline{d}^f_\mu(x)=0$ on a positive measure set. Indeed, the assumptions of Theorem~\ref{thprinc2} are not strong enough to ensure the almost everywhere existence of the pointwise dimension for the observations. The following result guaranties the existence for a large class of systems.
\begin{theorem}\label{thfcinf}
Let $f:\R^M\rightarrow\R^N$ be a $C^\infty$ function, let $\mu$ be any absolutely  continuous measure on $\R^M$. Then, $d^f_\mu$ exists and belongs to $\left\{0,1,...,\min\{M,N\}\right\}$ $\mu$-almost everywhere. More precisely, $d^f_\mu(x)=\rank d_x f$ for $\mu$-almost every $x\in\R^M$.
\end{theorem}
This is a non trivial result because the image measure $\pf$ may be quite complicated and rather counter intuitive. Already in the one dimensional case, there exists $f\in C^\infty(\R,\R)$ such that $f(\{f'=0\})$ is an uncountable set of dimension $0$ and $f_*(Leb|_{\{f'=0\}})$ is a non null and non atomic measure. We emphasize that Theorem~\ref{thfcinf} applies to any $C^\infty$ function, and not only for generic functions. This is essential in applications, where we are mostly interested in particular observables.
\begin{corollary}\label{cor2}
Let $T:X\subset\R^M\rightarrow X$ preserves an absolutely continuous invariant probability measure on $\R^M$ with super-polynomial decay of correlations. Let $f:\R^M\rightarrow\R^N$ be a Lipschitz $C^\infty$ observable. Then $R^f$ exists and belongs to $\left\{0,1,...,\min\{M,N\}\right\}$ almost everywhere.
\end{corollary}
\begin{proof}[Proof of Corollary~\ref{cor2}]
We apply Theorem~\ref{thprinc} and Theorem~\ref{thfcinf} when $d^f_\mu>0$. When $d^f_\mu=0$, we use Theorem~\ref{thsup}. $\Box$
\end{proof}


\subsection{On the necessity of the non-instantaneous recurrence rate}
In this part, we give a simple example which illustrates the utility of non-instantaneous return times.

Let $\Omega:=\{0,1\}^\N$ and $\sigma$ be the shift on $\Omega$. Fix some 1-approximable $\alpha\in\R$ (e.g. \cite{MR2112110} for a nice perspective) i.e. $\delta(\alpha)=1$ where
\[\delta(\alpha)=\sup\left\{\delta\geq 1,\,|\alpha-\frac{p}{q}|<\frac{1}{q^{1+\delta}}\textrm{ for infinitely many $\frac{p}{q}\in\Q$}\right\}.\]
Let $\nu$ be an invariant ergodic probability measure on $\Omega$. Fix some measurable $A\subset\Omega$ such that $1>\nu(A)>0$ and set $\varphi$:
\begin{equation}
\varphi(\omega)=\left\{ \begin{array}{ll}
0 & \textrm{ if $\omega\notin A$}\\
\alpha & \textrm{ if $\omega\in A$}.
\end{array}\right.
\end{equation}
Let $\T^1$ denote the 1-dimensional torus and define on $X:=\Omega\times\T^1$ the map 
\begin{eqnarray*}
T:&X&\longrightarrow X\\
&(\omega,y)&\longrightarrow (\sigma \omega, y+\varphi(\omega)).
\end{eqnarray*}
Let $Leb$ be the Lebesgue measure on $\T^1$. We consider the $T$-invariant probability measure $\mu:=\nu\otimes Leb$. We examine below the recurrence rate of the system $(X,T,\mu)$ for the observable $f$ given by the projection on the second variable i.e.
\begin{eqnarray*}
f:&X&\rightarrow\T^1\\
&(\omega,y)&\rightarrow y.
\end{eqnarray*}
First, we need the following obvious result on the pushforward measure: since $\pf=Leb$ and the local dimension of the Lebesgue measure is one, the measure $\pf$ is exact dimensional and satisfies
\begin{equation}\label{dfmu}\forall x\in X,\,\,d^f_\mu(x)=1.\end{equation}
\begin{prop}We have $R^f_i\neq d^f_\mu$ on a set of positive measure. More precisely
\[\forall x=(\omega,y)\in \Omega\backslash A\times\T^1,\,\,R^f_i(x)=0.\]
\end{prop}
\begin{proof} Let $\omega\in \Omega\backslash A$ and $y\in\T^1$, we have
\begin{eqnarray*}
f\left(T(\omega,y)\right)&=&f\left(\sigma\omega,y\right)\qquad\textrm{because $\omega\notin A$}\\
&=&y\\
&=&f\left(\omega,y\right).
\end{eqnarray*}
So, for all $r>0$, $\tau_r^f(x)=1$ and then $R^f_i(x)=0$.
\end{proof}
We therefore need to introduce the non-instantaneous return time to avoid this kind of problem.
\begin{prop}We have $R^f= d^f_\mu$ on a set of full measure.\end{prop}
\begin{proof} For $k\in\N$ and $\omega\in\Omega$, let $q_k(\omega):=\sum_{i=0}^{k-1}\bold{1}_A(\sigma^i\omega)$. Let $\eps>0$ fixed.  For $x=(\omega,y)$ and $n\in\N$, we have
\begin{eqnarray}
\tau_{\frac{1}{n^{1+\eps}},p}^f(x)&=&\inf\left\{k>p,\,f(T^kx)\in B\left(f(x),\frac{1}{n^{1+\eps}}\right)\right\} \nonumber \\
&=&\inf\left\{k>p,\,y+\alpha q_k(\omega)\in B\left(y,\frac{1}{n^{1+\eps}}\right)\right\} \nonumber \\
&=&\inf\left\{k>p,\,\|\alpha q_k(\omega)\|\leq\frac{1}{n^{1+\eps}}\right\} \label{sombirk}
\end{eqnarray}
where for $q\in\Z$
\[\|q\alpha\|:=\min\left\{|q\alpha-p|\, : \,p\in\Z\right\}.\]
Thanks to the choice of $\alpha$, there exists $k_0\in \N$ such that for all $k\geq k_0$, we get $\|k\alpha\|\geq\frac{1}{k^{1+\eps}}$. Taking $n\geq k_0$ we have
\begin{equation}\label{mneps}
m_{n,\eps}:=\inf\left\{q>k_0\, ,\,\| q\alpha\|\leq\frac{1}{n^{1+\eps}}\right\}\geq n.
\end{equation}
Since $\nu$ is ergodic, the Poincar\'e Recurrence Theorem gives, for $\nu$-almost every $\omega\in\Omega$
\begin{equation}\label{birk}
q_k(\omega)\underset{k\rightarrow\infty}{\longrightarrow}+\infty.
\end{equation}
So, for $\nu$-almost every $\omega\in\Omega$, we can choose $p$ sufficiently large such that $p\geq k_0$ and $q_p(\omega)\geq k_0$. If $k\in\N$ satisfies $p\leq k<m_{n,\eps}$, then $k_0\leq q_p(\omega)\leq q_k(\omega)\leq k<m_{n,\eps}$ and so $\|\alpha q_k(\omega)\|>\frac{1}{n^{1+\eps}}$. Since this is true for every $k\in[p,...,m_{n,\eps}]$ we end up with $\tau_{\frac{1}{n^{1+\eps}},p}^f(x)\geq m_{n,\eps}\geq n$. Finally, if $p\geq m_{n,\eps}$, we obviously have $\tau_{\frac{1}{n^{1+\eps}},p}^f(x)\geq m_{n,\eps}\geq n$.
Thus for $\mu$-almost every $x=(\omega,y)\in X$, we have
\begin{eqnarray*}
\underline{R}^f(x)&=&\lim_{p\rightarrow\infty}\liminf_{n\rightarrow+\infty}\frac{\log\tau_{\frac{1}{n^{1+\eps}},p}^f(x)}{-\log \frac{1}{n^{1+\eps}}}\\
&\geq&\liminf_{n\rightarrow+\infty}\frac{\log m_{n,\eps}}{-\log \frac{1}{n^{1+\eps}}}\\
&\geq&\lim_{n\rightarrow\infty}\frac{\log n}{\log n^{1+\eps}}=\frac{1}{1+\eps}.
\end{eqnarray*}
This is true for all $\eps>0$, thus
\[\underline{R}^f(x)\geq 1.\]
The conclusion follows from Theorem~\ref{thsup} and equation (\ref{dfmu}).
\end{proof}
\begin{remark}We point out that indeed our example fulfills the conditions of Corollary~\ref{thprinc} when, for example, $\nu$ is a Gibbs measure \cite{MR1919377}.\end{remark}


\section{Majoration of the recurrence rate for measure preserving systems}\label{majoration}
The basic strategy of the proof of Theorem~\ref{thsup} follows \cite{MR1833809}. We recall that any probability measure on $\R^N$ is weakly diametrically regular \cite{MR1833809}:
\begin{defi}
A measure $\mu$ is weakly diametrically regular (wdr) on the set $Z\subset X$ if for any $\eta>1$, for $\mu$-almost every $x\in Z$ and every $\eps>0$, there exists $\delta>0$ such that if $r<\delta$ then $\mu\left(B\left(x,\eta r\right)\right)\leq\mu\left(B\left(x,r\right)\right)r^{-\eps}$.
\end{defi}
\begin{proof}[Proof of Theorem~\ref{thsup}] The measure $\pf$ is weakly diametrically regular on $\R^N$. We can remark that the function $\delta(f(\cdot),\eps)$ in the previous definition can be made measurable for every fixed $\eps$. Let us fix $\eps>0$ and choose $\delta>0$ sufficiently small to have $\mu(G)>\mu(X)-\eps=1-\eps$ where
\[G:=\left\{x\in X,\,\delta(f(x),\eps)>\delta\right\}.\]
For all $r>0$, $\lambda>0$, $p\in\N$ and $x\in X$ we define the set
\[A_{r,x}:=\left\{y\in X,\,f(y)\in B\left(f(x),4r\right):\,\tau_{4r,p}^f(y,x)\geq\lambda^{-1}\pf\left(B\left(f(x),4r\right)\right)^{-1}\right\}\]
where $\tau^f_{4r,p}(y,x):=\inf\left\{k>p,\,d\left(f(T^ky),f(x)\right)<4r\right\}$ for $y\in f^{-1}B(f(x),4r)$. Markov's inequality gives:
\begin{equation}\label{cheb}
\mu(A_{r,x})\leq \lambda\pf\left(B\left(f(x),4r\right)\right)\int_{f^{-1}B(f(x),4r)}\tau^f_{4r,p}(y,x)\,d\mu(y).
\end{equation}
Since $\tau^f_{4r,p}(y,x)$ is bounded by the $p^{th}$ return time of $y$ in the set $f^{-1}B(f(x),4r)$, by Kac's lemma we have:
\begin{equation}\label{kacbosh}
\int_{f^{-1}B(f(x),4r)}\tau^f_{4r,p}(y,x)\,d\mu(y)\leq p.
\end{equation}
Using (\ref{cheb}) and (\ref{kacbosh}), we have:
\begin{equation}\label{ineg}
\mu(A_{r,x})\leq p\lambda\pf\left(B\left(f(x),4r\right)\right).\end{equation}
If $d(f(x),f(y))<2r$ then
\begin{equation}\label{inegarx}
\tau^f_{4r,p}(y,x)\pf\left(B\left(f(x),4r\right)\right)\geq\tau^f_{6r,p}(y)\pf\left(B\left(f(y),2r\right)\right).
\end{equation}
\begin{defi} Given $r>0$, a countable set $E\subset F$ is a maximal $r$-separated set for $F$  if
\begin{enumerate}
\item $B(x,\frac{r}{2})\cap B(y,\frac{r}{2})=\emptyset$ for any two distinct $x,y\in E$.
\item $\mu(F\,\backslash \underset{x\in E}{\bigcup}B(x,r))=0$.
\end{enumerate}
\end{defi}
Let $C\subset f(G)$ a maximal $2r$-separated set for $f(G)$.
\begin{eqnarray*}
D_\eps(r)&:=&\mu\left(\left\{y\in G, \tau^f_{6r,p}(y)\pf\left(B\left(f(y),2r\right)\right)\geq r^{-2\eps}\right\}\right)\\
&\leq&\sum_{f(x)\in C}\mu\left(\left\{y\in f^{-1}B\left(f(x),2r\right):\,\tau^f_{6r,p}(y)\pf\left(B\left(f(y),2r\right)\right)\geq r^{-2\eps}\right\}\right)\\
&\leq&\sum_{f(x)\in C}\mu(A_{r,x})\qquad\textrm{by (\ref{inegarx})}\\
&\leq&p\,r^{2\eps}\sum_{f(x)\in C}\pf\left(B(f(x),4r)\right)\qquad\textrm{with $\lambda=r^{2\eps}$ in (\ref{ineg})}\\
&\leq&p\,r^{\eps}\sum_{f(x)\in C}\pf\left(B(f(x),r)\right)\qquad\textrm{since $\pf$ is wdr and with $\eta=4$}\\
&\leq&p\,r^\eps\qquad\textrm{according to the definition of $C$}.
\end{eqnarray*}
Finally:
\[\sum_{n,e^{-n}<\delta}D_\eps(e^{-n})=\sum_{n>-\log\delta}D_\eps(e^{-n})\leq p\sum_ne^{-\eps n}<\infty.\]
Then, thanks to the Borel-Cantelli lemma, for $\mu$-almost every $x\in G$
\[\tau^f_{6e^{-n},p}(x)\pf\left(B\left(f(x),2e^{-n}\right)\right)\leq e^{2\eps n}\]
for any $n$ sufficiently large. Then
\begin{equation}\label{inegfin}
\frac{\log\tau^f_{6e^{-n},p}(x)}{n}\leq 2\eps+\frac{\log\pf(B(f(x),2e^{-n}))}{-n}.
\end{equation}
Observing that for all $a>0$ we have:
\begin{eqnarray*}
\underline{d}^f_\mu(x)=\underset{n\rightarrow\infty}{\underline\lim}\frac{\log\pf\left(B\left(f(x),ae^{-n}\right)\right)}{-n}&\textrm{ and }&\overline{d}^f_\mu(x)=\underset{n\rightarrow\infty}{\overline\lim}\frac{\log\pf\left(B\left(f(x),ae^{-n}\right)\right)}{-n}\\
\underline{R}^f(x)=\lim_{p\rightarrow\infty}\liminf_{n\rightarrow\infty}\frac{\log\tau_{ae^{-n},p}^f(x)}{n}&\textrm{ and }&\overline{R}^f(x)=\lim_{p\rightarrow\infty}\limsup_{n\rightarrow\infty}\frac{\log\tau_{ae^{-n},p}^f(x)}{n}
\end{eqnarray*}
and since $\eps$ can be chosen arbitrarily small , we have the result if we take the limit inferior or the limit superior and then the limit over $p$ in (\ref{inegfin}). \end{proof}
\section{Recurrence rate and dimension for mixing systems}\label{egalite}
Despite some similarities with \cite{MR2191396}, we emphasize that the proof of Theorem~\ref{thprinc2} is relatively different. In particular we make no assumption on the entropy of the system.
%
%
\begin{lem}\label{lemd>0}Under the assumptions of Theorem~\ref{thprinc2}, \[\textrm{$\underline{R}^f(x)>0$ for $\mu$-almost every $x$ such that $\underline{d}^f_\mu(x)>0$.}\]
\end{lem}
\begin{proof}
Let $X_+:=\{\underline{d}^f_\mu>0\}$. Let $\eps>0$ and let $a>0$ such that $\mu(X^+)\geq\mu(X_a)>\mu(X^+)-\eps$ where $X_a:=\{x\in X\,,\,\underline{d}^f_\mu(x)>a\}$.

We fix $b>0$ and for $\rho>0$ we consider the set $G=G_1\cap G_2\cap G_3$ with:
\begin{eqnarray*}
G_1&=&\left\{x\in X_a,\,\forall r\leq \rho,\,f_*\mu\left(B\left((f(x),2r\right)\right)\leq r^a\right\}\\
G_2&=&\left\{x\in X_a,\,\forall r\leq \rho,\,\pf\left(B\left((f(x),\frac{r}{2}\right)\right)\geq r^{N+b}\right\}\\
G_3&=&\left\{x\in X_a,\,\forall r\leq \rho,\,\pf\left(B\left((f(x),\frac{r}{2}\right)\right)\geq \pf\left(B\left(f(x),4r\right)\right)r^{a/2}\right\}.
\end{eqnarray*}
We observe that 
\begin{equation}\label{mug2}
\mu(G)\underset{\rho\rightarrow0}{\longrightarrow}\mu(X_a).
\end{equation}
Indeed, by definition of $\underline{d}^f_\mu$, we have $\mu(G_1)\rightarrow\mu(X_a)$. Moreover, since $\overline{d}^f_\mu\leq N$, $\mu(G_2)\rightarrow \mu(X_a)$, and since the measure $\pf$ is weakly diametrically regular, $\mu(G_3)\rightarrow\mu(X_a)$.
Let us define, for $n\in\N^*$, $\eps_n:=\frac{1}{n^{4/a}}$ and
\[A_n:=\left\{x\in X\,,\,f(T^nx)\in B\left(f(x),\eps_n\right)\right\}.\]
Let $x\in G$ and $n\in\N^*$, we have
\[f^{-1}B\left(f(x),\eps_n\right)\cap A_n\subset f^{-1}B\left(f(x),\eps_n\right)\cap T^{-n}f^{-1}B\left(f(x),2\eps_n\right).\]
Let $\eta_{\eps_n}:[0,+\infty)\rightarrow\R$ be the $\frac{1}{\eps_n}$-Lipschitz map such that $1_{[0,\eps_n]}\leq\eta_{\eps_n}\leq1_{[0,2\eps_n]}$ and set $\psi_{x,\eps_n}(y)=\eta_{\eps_n}\left(d(f(x),f(y))\right)$.
 Since $f$ is $L$-Lipschitz, $\psi$ is clearly $\frac{L}{\eps_n}$-Lipschitz. Since $(X,T,\mu)$ has super-polynomial decay of correlation, we have:
\begin{eqnarray*}
& &\mu\left(f^{-1}B(f(x),\eps_n)\cap T^{-n}f^{-1}B(f(x),2\eps_n)\right)\leq\int_X\pep\pdep d\mu(y)\\
& &\leq \|\psi_{x,\eps_n}\|\|\psi_{x,2\eps_n}\|\theta_n+\int_X\pep d\mu(y)\int_X\pdepy d\mu(y)\\
& &\leq\frac{L^2}{\eps_n^2}\theta_n+\pf(B(f(x),2\eps_n))\pf(B(f(x),4\eps_n)).
\end{eqnarray*}
There exists $n_0\in\N$ such that $\forall n\geq n_0$, we have $\eps_n=\frac{1}{n^{4/a}}<\rho$ and using the definition of $G$, for all $n\geq n_0$:
\begin{eqnarray*}
\mu\left(f^{-1}B\left(f(x),\eps_n\right)\cap A_n\right)&\leq&L^2\eps_n^{-2-N-b}\theta_n\pf\left(B\left(f(x),\frac{\eps_n}{2}\right)\right)+\\& &\qquad+\eps_n^{a-a/2}\pf\left(B\left(f(x),\frac{\eps_n}{2}\right)\right)\\
&\leq&\pf\left(B\left(f(x),\frac{\eps_n}{2}\right)\right)\left[L^2(\eps_n)^{-2-N-b}\theta_n+\eps_n^{a/2}\right].
\end{eqnarray*}
Let $B\subset G$ such that $\left(B\left(f(x),\eps_n\right)\right)_{x\in B}$ is a maximal $\eps_n$-separated set for $f(G)$. Since $\left(f^{-1}B\left(f(x),\eps_n\right)\right)_{x\in B}$ covers $G$, we have:
\begin{eqnarray*}
\mu\left(G\cap A_n\right)&\leq&\sum_{x\in B}\pf\left(B\left(f(x),\eps_n\right)\cap A_n\right)\\
&\leq&\sum_{x\in B}\pf\left(B\left(f(x),\frac{\eps_n}{2}\right)\right)\left[L^2(\eps_n)^{-2-N-b}\theta_n+
\eps_n^{a/2}\right]\\
&\leq&L^2(\eps_n)^{-2-N-b}\theta_n+\eps_n^{a/2}.
\end{eqnarray*}
Since $\sum_{n\in\N^*}\eps_n^{a/2}=\sum_{n\in\N^*}\frac{1}{n^{2}}<\infty$ and since the decay of correlations is super-polynomial, we obtain:
\begin{equation}
\sum_{n\in \N^*}\mu\left(G\cap A_n\right)<+\infty.
\end{equation}
By Borel-Cantelli lemma and using (\ref{mug2}), we have that for $\mu$-almost every $x\in X_a$, there exists $n_1(x)$ such that for every $n\geq n_1(x)$, $f(T^nx)\notin B\left(f(x),\frac{1}{n^{4/a}}\right)$. So, for  $\mu$-almost every $x\in X_a$, for $p\geq n_1(x)$ and $n\geq n_1(x)$,
\begin{equation}
\tau_{\frac{1}{n^{4/a}},p}^f(x)>n
\end{equation}
which gives us
\begin{eqnarray*}
\underline{R}^f(x)&=&\lim_{p\rightarrow\infty}\liminf_{r\rightarrow0}\frac{\log\tau_{r,p}^f(x)}{-\log r}\\
&=&\lim_{p\rightarrow\infty}\liminf_{n\rightarrow+\infty}\frac{\log\tau_{\frac{1}{n^{4/a}},p}^f(x)}{-\log \frac{1}{n^{4/a}}}\\
&\geq&\lim_{n\rightarrow+\infty}\frac{\log n}{\log n^{4/a}}=\frac{a}{4}>0.
\end{eqnarray*}
Since we can choose $\eps$ arbitrarily small, the lemma is proved. \end{proof}

%
%
%
\begin{lem}\label{lemintervalle}
Let $a>0$, $\delta>0$ and $\eps>0$. Let $X_a:=\{x\in X\,,\,\underline{d}^f_\mu(x)>a\}$. For $\mu$-almost every $x\in X_a$, there exists $r(x)>0$ such that for every $r\in]0,r(x)[$ and for every integer $n\in[r^{-\delta},\pf\left(B(f(x),r)\right)^{-1+\eps}]$, we have $d\left(f(T^nx),f(x)\right)\geq r$.
\end{lem}
\begin{proof} Let $\delta>0$ and $\eps>0$. We fix $b>0$, $c=\frac{a\eps}{2}$ and for $\rho>0$ we consider the set $G=G_1\cap G_2\cap G_3$ with:
\begin{eqnarray*}
G_1&=&\left\{x\in X_a,\,\forall r\leq \rho,\,f_*\mu\left(B\left((f(x),2r\right)\right)\leq r^a\right\}\\
G_2&=&\left\{x\in X_a,\,\forall r\leq \rho,\,\pf\left(B\left((f(x),\frac{r}{2}\right)\right)\geq r^{N+b}\right\}\\
G_3&=&\left\{x\in X_a,\,\forall r\leq \rho,\,\pf\left(B\left((f(x),\frac{r}{2}\right)\right)\geq \pf\left(B\left(f(x),4r\right)\right)r^c\right\}.
\end{eqnarray*}
While proving Lemma~\ref{lemd>0}, we already observed that 
\begin{equation*}
\mu(G)\underset{\rho\rightarrow0}{\longrightarrow}\mu(X_a).
\end{equation*}
Let $r\leq \rho$, we define:
\[A_\eps(r):=\left\{y\in X\,:\,\exists n\in[r^{-\delta},\pf\left(B(f(y),r)\right)^{-1+\eps}]\textrm{ such that }d\left(f(T^ny),f(y)\right)<r\right\}.\]
Let $x\in G$, we have:
\begin{eqnarray*}
& &f^{-1}B(f(x),r)\cap A_\eps(r)\\
&=&\left\{y,f(y)\in B(f(x),r),\exists n\in[r^{-\delta},\pf\left(B(f(y),3r)\right)^{-1+\eps}],d(f(T^ny),f(y))<r\right\}\\
&\subset&\left\{y,f(y)\in B(f(x),r),\exists n\in[r^{-\delta},\pf\left(B(f(x),2r)\right)^{-1+\eps}],d(f(T^ny),f(x))<2r\right\}\\
&=&\underset{r^{-\delta}\leq n\leq \pf\left(B(f(x),2r)\right)^{-1+\eps}}{\bigcup}f^{-1}B(f(x),r)\cap T^{-n}f^{-1}B(f(x),2r).
\end{eqnarray*}
Let $\eta_{r}:[0,+\infty)\rightarrow\R$ be the $\frac{1}{r}$-Lipschitz map such that $1_{[0,r]}\leq\eta_r\leq1_{[0,2r]}$ and set $\psi_{x,r}(y)=\eta_{r}\left(d(f(x),f(y))\right)$.
 Since $f$ is $L$-Lipschitz, $\psi$ is clearly $\frac{L}{r}$-Lipschitz. Using the assumption on the decay of correlations of $(X,T,\mu)$, we obtain
\begin{eqnarray*}
& &\mu\left(f^{-1}B(f(x),r)\cap T^{-n}f^{-1}B(f(x),2r)\right)\leq\int_X\psi_{x,r}(y)\psi_{x,2r}(T^ny) d\mu(y)\\
& &\leq \|\psi_{x,r}\|\|\psi_{x,2r}\|\theta_n+\int_X\psi_{x,r}(y) d\mu(y)\int_X\psi_{x,2r}(y) d\mu(y)\\
& &\leq\frac{L^2}{r^2}\theta_n+\pf(B(f(x),2r))\pf(B(f(x),4r)).
\end{eqnarray*}
Let us choose $k>1$ such that $\delta(k-1)-2\geq N+2b$ and we choose $\rho$ such that $n\geq\rho^{-\delta}$ implies $(k-1)(n+1)^{-k}\geq\theta_n$ (which is possible by definition of $\theta_n$). Let $r\in(0,\rho)$, set $I_r=[r^{-\delta},\pf\left(B(f(x),2r)\right)^{-1+\eps}]\cap\N$, we have
\begin{eqnarray*}
\mu\left(f^{-1}B(f(x),r)\cap A_\eps(r)\right)&\leq&\underset{n\in I_r}{\sum}\frac{L^2}{r^2}\theta_n+\pf(B(f(x),2r))\pf(B(f(x),4r))\\
&\leq&\frac{r^{\delta(k-1)-2}}{L^2}+\pf(B(f(x),2r))^\eps\pf(B(f(x),4r))\\
&\leq& \frac{r^{N+2b}}{L^2}+r^{a\eps}\pf(B(f(x),\frac{r}{2}))r^{-c}\qquad\textrm{by definition of $G$}\\
&\leq& \pf(B(f(x),\frac{r}{2}))\left(\frac{r^b}{L^2}+r^{a\eps-c}\right).
\end{eqnarray*}
Let $B\subset G$ such that $\left(f(x)\right)_{x\in B}$ is a maximal $r$-separated set for $f(G)$. Since the collection $\left(f^{-1}B\left(f(x),r\right)\right)_{x\in B}$ covers $G$, we have:
\begin{eqnarray*}
\mu\left(G\cap A_\eps(r)\right)&\leq&\sum_{x\in B}\mu\left(f^{-1}B\left(f(x),r\right)\cap A_\eps(r)\right)\\
&\leq&\sum_{x\in B}\pf\left(B\left(f(x),\frac{r}{2}\right)\right)\left(\frac{r^b}{L^2}+r^{a\eps-c}\right)\\
&\leq&\frac{r^b}{L^2}+r^{a\eps/2}\qquad\textrm{since $B\left(f(x),\frac{r}{2}\right)$ are disjoints}.
\end{eqnarray*}
Then 
\[\sum_{k\in\N}\mu\left(G\cap A_\eps(e^{-k})\right)<+\infty\]
thus, by Borel-Cantelli lemma, we have for $\mu$-almost every $y\in G$, there exists $n_1(y)$ such that for every $k\geq n_1(y)$, $y\notin A_\eps(e^{-k})$. So, for $r$ sufficiently small there exists $k\in\N$ such that $e^{-k-1}<r\leq e^{-k}\leq e^{-n_1(y)}$ and since $e^{\delta k}\leq r^{-\delta}$ and $3e^{-m}<3er$, there does not exist $n\in[r^{-\delta},\pf\left(B(f(y),3r)\right)^{-1+\eps}]$ such that $d(f(T^ny),f(y))<r$. Since $\pf$ is weakly diametrically regular the factor $3e$ is irrelevant and the lemma is proved. \end{proof}
%
%
%
%
\begin{proof}[Proof of Theorem~\ref{thprinc2}]Let $\zeta>0$, since $ \underline{R}^f(x)>0$ for $\mu$-almost every $x\in X^+:=\{\underline{d}^f_\mu>0\}$ by Lemma~\ref{lemd>0}, there exists $a>0$ such that $\mu(X^+)\geq\mu(\{ \underline{R}^f>a\})>\mu(X^+)-\zeta$. For any $x\in\{ \underline{R}^f>a\}$ , for $p$ sufficiently large and $r$ sufficiently small, we have 
\[\tau_{r,p}^f(x)\geq r^{-a}.\]
Thanks to Lemma~\ref{lemintervalle} with $\delta=a$ and $\eps>0$, for $\mu$-almost every $x\in \{ \underline{R}^f>a\}$, if $r$ is sufficiently small and $p$ sufficiently large, then $\tau_{r,p}^f(x)\geq f_*\mu\left(B(f(x),r)\right)^{-1+\eps}$. Thus, $\underline{R}^f\geq(1-\eps)\underline{d}^f_\mu$ and $\overline{R}^f\geq(1-\eps)\overline{d}^f_\mu $ $\mu$-almost everywhere on $\{\underline{R}^f>a\}$. The theorem is proved choosing $\eps>0$ arbitrarily small and then $\zeta>0$ arbitrarily small. \end{proof}
%
%
%
%
\section{Dimensions of the smooth image of Lebesgue measure}\label{dimensions}
\subsection{Hausdorff and packing dimensions}\label{sechaus}
In this section, we recall the notion of Hausdorff dimension, packing dimension and pointwise dimension and the link between each other (see \cite{MR1449135} for more details).

Let $(X,d)$ be a metric space. Let $U$ be a non-empty set, its diameter is
\[\textrm{diam }U:=\sup\left\{d(x,y)\,:\,x,y\in U\right\}.\]
Given $\delta>0$, a collection $\{U_i\}_{i\in I}$ is called a countable $\delta$-cover of a set $E$ if $I$ is countable, $E\subset\cup_{i\in I}U_i$ and for all $i\in I$, $0<\textrm{diam }U_i\leq\delta$.

Let $E$ be a subset of $X$ and $s\geq0$, for $\delta>0$, we define:
\begin{equation}
\mathcal{H}^s_\delta(E)=\inf\left\{\sum_{i\in I}(\textrm{diam }U_i)^s\,:\,\{U_i\}_{i\in I}\textrm{ is a countable $\delta$-cover of $E$}\right\}.
\end{equation}
We then define the \textit{Hausdorff $s$-dimensional outer measure of $E$} as
\begin{equation}
\mathcal{H}^s(E)=\lim_{\delta\rightarrow0}\mathcal{H}^s_\delta(E).
\end{equation}
There exists a unique $t$ such that $\mathcal{H}^s(E)=\infty$ if $s<t$ and $\mathcal{H}^s(E)=0$ if $s>t$ which is called the \textit{Hausdorff dimension of $E$} i.e.
\begin{equation}
\dim_HE=\inf\left\{s\,:\,\mathcal{H}^s(E)<\infty\right\}=\sup\left\{s\,:\,\mathcal{H}^s(E)>0\right\}.
\end{equation}
If $\mu$ is a probability measure on $X$, we define the \textit{Hausdorff dimension of $\mu$}
\begin{equation}
\dim_H\mu=\inf \left\{\dim_HY\,:\,\mu(Y^c)=0\right\}.
\end{equation}
\begin{remark}We warn the reader that this definition of the Hausdorff dimension of a measure differs from the one given by Falconer \cite{MR1449135} but it is the most used in Ergodic Theory.\end{remark}

Given $\eps>0$, the collection $\left\{B(x_i,r_i)\right\}_{i\in I}$ is called a $\eps$-packing of $E$ if $I$ is a finite or countable set, for all $i\in I$ we have $x_i\in E$, $r_i\leq\eps$ and the balls are disjoints. For $s\geq0$, we write
\[\mathcal{P}^s_\eps(E)=\sup\left\{\sum_{i\in I}(r_i)^s\,:\,\{B(x_i,r_i)\}_{i\in I}\textrm{ is a $\eps$-packing of $E$}\right\}\]
and 
\[\mathcal{P}^s_0(E)=\lim_{\eps\rightarrow0}\mathcal{P}^s_\eps(E).\]
We then introduce the \textit{$s$-dimensional packing outer measure}
\begin{equation}\label{sdimpacmeas}
\mathcal{P}^s(E)=\inf\left\{\sum_{i=1}^\infty\mathcal{P}_0^s(E_i)\,:\,E\subset\bigcup_{i=1}^\infty E_i\right\}
\end{equation}
and the \textit{packing dimension of $E$} is defined as Hausdorff dimension
\begin{equation}
\dim_PE=\inf\left\{s\,:\,\mathcal{P}^s(E)<\infty\right\}=\sup\left\{s\,:\,\mathcal{P}^s(E)>0\right\}.
\end{equation}
For a probability measure $\mu$, we also have a \textit{packing dimension of $\mu$}
\begin{equation}
\dim_P\mu=\inf \left\{\dim_PY\,:\,\mu(Y^c)=0\right\}.
\end{equation}
There is a link between Hausdorff dimension, packing dimension and pointwise dimension:
\begin{prop}Assume that $X\subset\R^N$ for some $N$,
\begin{equation}\label{essup}\dim_H\mu=\textrm{ess-sup }\underline{d}_\mu\end{equation}
and
\begin{equation}\label{dimpacesssup}\dim_P\mu=\textrm{ess-sup } \overline{d}_\mu.\end{equation}
\end{prop}
%
%
%
%
\subsection{Existence of the pointwise dimension}
Bates and Moreira proved a generalization of the classical Morse-Sard Theorem for Hausdorff measures. Unfortunately, in view of \eqref{essup}, this is not enough to get an upper bound for the upper pointwise dimension. A key ingredient of their proof is the following generalized Morse decomposition. Given a differentiable $f$ from $\R^M$ to $\R^N$, for $\kappa\in\{0,1,...,\min\{M,N\}\}$, we define $C_\kappa:=\{x\in\R^M\,,\, \textrm{rank}(d_xf)=\kappa\}$.
\begin{lem}[\cite{MR1805620}]\label{decomposition}
Let $f\in C^k(\R^M,\R^N)$ with $k\geq 2$. Let $\kappa\leq M$ be an integer. Let $\eta>0$. There is a decomposition $\{A_i\}_{i\in\N}$ of $C_\kappa$ such that for each $i\in \N$ there exist two subspaces $E'_i$ and $E''_i$ which satisfy $\R^M=E'_i\oplus E''_i$, $\dim E'_i\leq \kappa$ and if $S\subset\R^M$:
\begin{equation}\label{eqdiam}
\diam\left(f(S\cap A_i)\right)\leq \left(\|f|_{A_i}\|_{C^1}+\eta\right)\diam(\pi_{E'_i}S)+\eta(\diam S)^k.
\end{equation}
\end{lem}
This decomposition will be instrumental to prove an analogue result but for the packing dimension.
\begin{lem}\label{lempack}
If $f\in C^\infty(\R^M,\R^N)$ then the packing dimension of the critical set satisfies
\[\dim_Pf(C_\kappa)\leq \kappa.\]
\end{lem}
\begin{proof}
Let $k\geq2$, since $f$ is of class $C^\infty$, $f$ is of class $C^k$. Let $\{A_i\}_{i\in\N}$ be the decomposition of $C_\kappa$ given by Lemma~\ref{decomposition} with $\eta=1$. Let $K\subset\R^M$ be a compact set. Let $i\in\N$. Let $d_i$ be the distance in $\R^M$ such that for $x\in \R^M$ and $y\in\R^M$, $d_i(x,y)=d(\pi_{E'_i}x,\pi_{E'_i}y)+d(\pi_{E''_i}x,\pi_{E''_i}y)$. Let $\eps>0$ and $\left\{B(f(x_j),r_j)\right\}_{j\in J}$ a $\eps$-packing of $f(A_i\cap K)$. Let $j\in J$ and $l\in J$, in (\ref{eqdiam}), we take $S:=\{x_j,x_l\}$, then
\begin{eqnarray*}
r_j+r_l&\leq& d_i(f(x_j),f(x_l))=\textrm{diam}(f(S\cap A_i))\\
&\leq& Cd(\pi_{E'_i}x_j,\pi_{E'_i}x_l)+ \left(d(\pi_{E'_i}x_j,\pi_{E'_i}x_l)+d(\pi_{E''_i}x_j,\pi_{E''_i}x_l)\right)^k
\end{eqnarray*}
where $C:=\|f\|_{C^1}+1$. This implies that
\begin{itemize}
\item either $Cd(\pi_{E'_i}x_j,\pi_{E'_i}x_l)\geq\frac{1}{2}(r_j+r_l)$
\item or $Cd(\pi_{E'_i}x_j,\pi_{E'_i}x_l)\leq\frac{1}{2}(r_j+r_l)$ and then $\left(d(\pi_{E'_i}x_j,\pi_{E'_i}x_l)+d(\pi_{E''_i}x_j,\pi_{E''_i}x_l)\right)^k\geq\frac{1}{2}(r_j+r_l)$. If $\eps$ is sufficiently small (depending only on $k$), we have
\[\frac{1}{2C}(r_j+r_l)\leq\left(1-\left(\frac{1}{2}\right)^{1/k}\right)\left(\frac{1}{2}(r_j+r_l)\right)^{1/k}\]
thus
\[d(\pi_{E''_i}x_j,\pi_{E''_i}x_l)\geq\left(\frac{1}{2}\right)^{1/k}\left(\frac{1}{2}(r_j+r_l)\right)^{1/k}.\]
\end{itemize}
For $j\in J$, let $S_j\subset E'_i\times E''_i$ be the product of the ball $B_{E_i'}(\pi_{E_i'}x_j,\frac{1}{4C}r_j)$ with the ball $B_{E_i''}(\pi_{E_i''}x_j,\frac{1}{2}\left(\frac{1}{4}r_j\right)^{1/k})$. If $l\neq  j$, we have $S_j\cap S_l=\emptyset$ since
\begin{itemize}
\item either $d(\pi_{E'_i}x_j,\pi_{E'_i}x_l)\geq\frac{1}{2C}(r_j+r_l)>\frac{1}{4C}r_j+\frac{1}{4C}r_l$
\item or $d(\pi_{E''_i}x_j,\pi_{E''_i}x_l)\geq\left(\frac{1}{4}(r_j+r_l)\right)^{1/k}>\frac{1}{2}\left(\frac{1}{4}r_j\right)^{1/k}+\frac{1}{2}\left(\frac{1}{4}r_l\right)^{1/k}$.
\end{itemize}
There exists a constant $\delta_i$ such that 
\[\diam S_j\leq\delta_i\eps^{1/k}.\]
The rectangles $S_j$ are disjoints and have non empty intersection with $K$, thus
\begin{equation}\label{volume}\sum_{j\in J}\textrm{Vol}(S_j)\leq \textrm{Vol}(K+B(0,\delta_i\eps^{1/k})).\end{equation}
Let $p=\dim E_i'$. There exists a constant $\gamma_i$ such that the volume of each $S_j$ is
\[\textrm{Vol}(S_j)=\gamma_i\left(\frac{1}{4C}r_j\right)^p\times\left[\frac{1}{2}\left(\frac{1}{4}r_j\right)^{1/k}\right]^{M-p}.\]
This implies together with \eqref{volume}
\[\sum_{j\in J}(r_j)^{p+(M-p)/k}\leq c(i,k,K)<\infty\]
where $c(i,k)$ is a finite constant depending on $i$,$k$ and $K$.\\
Now, by definition of  the $s$-dimensional packing measure (\ref{sdimpacmeas}), computed with the particular metric $d_i$, we obtain:
\[\mathcal{P}^{M/k+p(1-1/k)}(f(A_i\cap K))\leq c(i,k,K)<\infty.\]
This inequality holds for the packing measure computed with the metric $d_i$ and thus this is also true (possibly with another constant) for the packing measure computed with the euclidean metric $d$ since they are equivalent. Therefore
\[\dim_Pf(A_i\cap K)\leq\frac{M}{k}+p(1-\frac{1}{k})\leq \frac{M}{k}+\kappa.\]
Finally, taking a sequence of compacts $K_n$ such that $\R^M=\cup_{n\in\N}K_n$, we obtain:
\begin{eqnarray*}
\dim_Pf(C_\kappa)&=&\dim_Pf(\underset{i,n\in\N}{\bigcup} A_i\cap K_n)\\
&=&\dim_P\underset{i,n}{\bigcup}f(A_i\cap K_n)\\
&=&\sup_{i,n}\dim_Pf(A_i\cap K_n)\qquad\textrm{see \cite{MR1449135}}\\
&\leq&\frac{M}{k}+\kappa.
\end{eqnarray*}
Since $k$ is arbitrarily large we get
\begin{equation}\label{dimpaccr}
\dim_Pf(C_\kappa)\leq \kappa.
\end{equation}
\end{proof}
Without loss of generality, we prove Theorem~\ref{thfcinf} on $\R^M$ with $\mu$ equal to the Lebesgue measure $\lambda$. The general case can be deduced from it easily.
\begin{proof}[Proof of Theorem~\ref{thfcinf}] Let $\kappa\in\{0,...,\min\{M,N\}\}$.\\
\underline{1. If $A\subset C_\kappa$ and $\lambda(A)>0$ then $\dim_H f_*(\lambda|_A)\geq \kappa$:}\\
Indeed, let $B\subset C_\kappa$ with $\lambda(B)>0$, there exist $V$ open with $\lambda(B\cap V)>0$ and $\pi:\R^N\rightarrow\R^\kappa$ a linear map such that $d_x \pi\circ f$ is of maximal rank $\kappa$ for every $x\in V$ and so  $f_\kappa:=\pi \circ f$ satisfies $Jf_\kappa\neq0$ for every $x\in V$ (where $J$ is the Jacobian i.e. $Jf_\kappa=\sqrt{\det(d_xf_\kappa)(d_xf_\kappa)^t}$). Since $\pi$ is Lipschitz, it is known \cite{MR1449135} that:
\begin{equation}\label{dimproj}
\dim_H\pi(f(B))\leq\dim_Hf(B).
\end{equation}
Using the coarea formula (e.g. \cite{MR1158660}, in fact we could have worked directly with $f$ using \cite{MR528671}):
\begin{equation}\label{coarea}
\int_B Jf_\kappa d\lambda=\int_{f_\kappa(B)}H^{M-\kappa}\left( B\cap f_\kappa^{-1}(\{y\})\right)d\lambda_\kappa(y)
\end{equation}
where $\lambda_\kappa$ is the Lebesgue measure on $\R^\kappa$. Since $\lambda(B\cap V)>0$ and $Jf_\kappa(x)\neq0$ for every $x\in B\cap V$, the left-hand side of (\ref{coarea}) does not vanish and therefore, neither does the right-hand side.  Then $\lambda_\kappa(f_\kappa(B))>0$ and so $\dim_Hf_\kappa(B)\geq \kappa$ which gives, using (\ref{dimproj}):
\begin{equation}\label{dimhb}
\kappa\leq \dim_Hf(B).\end{equation}
Let $A\subset C_\kappa$ with $\lambda(A)>0$. We recall
\begin{equation}\label{defdimhf}
\dim_Hf_*(\lambda|_A):=\inf\{\dim_H Y\textrm{ such that }(f_*(\lambda|_A)(Y^c)=0\}.
\end{equation}
Let $Y$ be such that $(f_*(\lambda|_A)(Y^c)=0$, since $Y\supset f\left(A\cap f^{-1}(Y)\right)$ we have
\[\dim_HY\geq\dim_Hf\left(A\cap f^{-1}(Y)\right).\]
Moreover, since $\lambda(f^{-1}(Y)\cap A)=\lambda(A)>0$ and $f^{-1}(Y)\cap A\subset C_\kappa$, we can choose $B=f^{-1}(Y)\cap A$ in the previous consideration and (\ref{dimhb}) gives:
\[\dim_HY\geq\dim_Hf\left(A\cap f^{-1}(Y)\right)\geq \kappa\]
and then
\begin{equation}\label{dimhf}
\dim_Hf_*(\lambda|_A)\geq \kappa.
\end{equation}
We define  $\nu:=f_*\lambda$ and $\nu_\kappa:=f_*(\lambda|_{C_\kappa})$.\\
\underline{2. Let us prove that $d_\nu=\kappa$ $\nu_\kappa$-almost everywhere:}\\
$\bullet$ Firstly, since $\nu_\kappa$ is supported by $f(C_\kappa)$,  by Lemma~\ref{lempack} we have 
\[\dim_P\nu_\kappa\leq\dim_Pf(C_\kappa)\leq \kappa.\]
Since the packing dimension satisfies the relation (\ref{dimpacesssup}) we get
\[\overline{d}_{\nu_\kappa}(x)\leq \kappa\textrm{ for $\nu_\kappa$-almost every $x\in \R^N$}.\]
Since, for every $x\in \R^N$ and every $\eps>0$, $\nu\left(B(x,\eps)\right)\geq\nu_\kappa\left(B(x,\eps)\right)$, we have
\[\overline{d}_{\nu}(x)\leq\overline{d}_{\nu_\kappa}(x)\]
and then 
\begin{equation}\label{inegd1}\overline{d}_{\nu}(x)\leq \kappa\textrm{ for $\nu_\kappa$-almost every $x\in \R^N$}.\end{equation}
$\bullet$ Let $K\subset\R^N$ be a compact subset. Let $Z:=\left\{\underline{d}_\nu\leq\rho\right\}\cap K$ with $\rho<\kappa$. If $\nu_\kappa(Z)>0$ then $\lambda(C_\kappa\cap f^{-1}(Z))>0$, thus by (\ref{dimhf}) we obtain $\dim_Hf_*(\lambda|_{C_\kappa\cap f^{-1}(Z)})\geq \kappa$. By (\ref{defdimhf}) and since $f_*(\lambda|_{C_\kappa\cap f^{-1}(Z)})(Z^c)=0$, we obtain $\dim_HZ\geq \kappa$. 

On the other hand, by definition of $\underline{d}_\nu$, $\forall x\in Z$, $\exists J_x\subset\R^+$ with $0\in\overline{J_x}$, such that $\forall r\in J_x$, $\nu\left(B(x,r)\right)\geq r^\rho$. We notice that $\left\{B(x,r)\,,\,x\in Z\,,\,r\in J_x\cap[0,1]\right\}$ cover $Z$ so, by Besicovitch covering Theorem, there exists a subcovering $\left\{B(x_i,r_i)\right\}_{i\in I}$ with $I$ countable and $m_0$ a constant depending only on $N$ such that $Z\subset\cup_{i\in I}B(x_i,r_i)$ and the multiplicicty of the subcovering is bounded by $m_0$. Hence
\[\sum_{i\in I}r_i^\rho\leq\sum_{i\in I}\nu\left(B(x_i,r_i)\right)\leq m_0\nu(K)\]
which implies $\dim_HZ\leq\rho$. But this is in contradiction with the fact that $\dim_HZ\geq \kappa$. Then, for all $\rho<\kappa$ and for all compact $K$, $\nu_\kappa(\left\{\underline{d}_\nu\leq\rho\right\}\cap K)=0$. Thus
\begin{equation}\label{inegd2}\underline{d}_{\nu}(x)\geq \kappa\textrm{ for $\nu_\kappa$-almost every $x\in \R^N$}.\end{equation}
\underline{3. Conclusion:} Using (\ref{inegd1}) together with (\ref{inegd2}) implies that $d_\nu=\kappa$ $\nu_\kappa$-almost everywhere. The theorem follows from $\R^M=C_0\cup C_1\cup...\cup C_{\min\{M,N\}}$. \end{proof}

\bibliographystyle{siam} 
\bibliography{biblio}

\end{document}